\documentclass[11pt]{article}
\usepackage{tikz}
\usepackage[margin=1in]{geometry}
\usetikzlibrary{arrows}
\usepackage[T1]{fontenc}
\usepackage{textcomp}
\usepackage{color}
\usepackage{graphicx}

\usepackage[scaled=0.92]{helvet}
\usepackage{times}
\usepackage{amssymb}
\usepackage[all,knot,arc,curve,color,frame]{xy}
\usepackage{url}
\input xy
\xyoption{all}
\usepackage{extarrows}
\usepackage{enumitem}
\setlist{nosep}

\usepackage[all,knot,arc]{xy}
\usepackage{mathtools}

\newtheorem{theo}{Theorem}
\newtheorem{lemma}[theo]{Lemma}

\newtheorem{cor}[theo]{Corollary}

\title{Optimal Packings of 22 and 33 Unit Squares in a Square}
\author{Wolfram Bentz \\
\small Department of Physics and Mathematics \\
\small University of Hull \\
\small United Kingdom \\
\small {\tt W.Bentz@hull.ac.uk}} 
\begin{document}
\maketitle

\begin{abstract} 
Let $s(n)$ be the side length of the smallest square into which $n$ non-overlapping unit squares can be packed. In 2010, the author showed that $s(13)=4$ and $s(46)=7$. Together with the result $s(6)=3$ by Keaney and Shiu, these results strongly suggest that   
$s(m^2-3)=m$ for $m\ge 3$, in particular for the values $m=5,6$, which correspond to cases that lie in between the  previous results. 

In this article we show that indeed $s(m^2-3)=m$ for $m=5,6$,
implying that the most efficient packings of 22 and 33 squares are the trivial ones. To achieve our results, we modify the well-known method 
of sets of unavoidable points by replacing them with continuously varying families of such sets. 
\end{abstract}
\section{Introduction}

The study of packing unit squares 
into a square goes back to Erd\"os and Graham \cite{eg}, who 
examined the asymptotic packing efficiency as the side length of the containing square increased towards infinity. 
G\"obel \cite{goe} was the first to show that particular packings are optimal for a given non-square number of unit squares. The search for good packings for given number of unit squares 
was addressed  in the popular science literature in various articles by Gardner
\cite{mg}.   

Let $s(n)$ be the side length of the smallest square into which $n$ non-overlapping unit squares can be packed. Non-trivial cases
for which $s(n)$ is known are $s(m^2-1)=s(m^2-2)=m$ for $m \ge 2$ (Nagamochi 
\cite{naga}, single values previously shown by 
G\"obel \cite{goe}, El Moumni \cite{el}, and Friedman \cite{fried}), $s(5)=2+{\frac{1}{2}}\sqrt{2}$ (G\"obel \cite{goe}), $s(6)=3$ (Kearney and Shiu
 \cite{ks}), $s(10)=3+{\frac{1}{2}}\sqrt{2}$ 
(Stromquist \cite{strom}), $s(13)=4$, and $s(46)=7$ (Bentz \cite{b46}). There are moreover non-trivial best packings and lower bounds known for various values of $n$. 
Examples on many of these results and the underlying techniques used are given in the  survey article by Friedman \cite{fried}.

In \cite{ks} and \cite{b46}, it was shown that  $s(m^2-3)=m$ for $m=3,4,7$. These results suggested that the  holds for the intermediate values
 $m=5,6$. We will show this result in this article, by adopting the proof for $m=7$ from \cite{b46}.
 Previously, the best lower bounds in these cases are  $s(22) \ge \sqrt{15}+1\approx 4.87298$ and $s(33)\ge \sqrt{24}+1\approx 5.89898$, and follow  from a general result in Nagamochi \cite{naga}.

As the trivial (or ``chess board'') packings show that $s(22)\le 5$, and $s(33) \le 6$, it it suffices to establish the opposite inequality.  
Let 
a box be the interior of any square with side length $s$ satisfying $1<s\le 1.01$.  Following Stromberg, we will establish that $m^2-3$ squares
cannot be packed in a square with side length smaller then $m$ by proving the equivalent statement that it is impossible to pack $m^2-3$ boxes in a square of side length $m$.

In order to do so, we will adopt the previously used method of unavoidable points to continuously varying sets of such points.
We will introduce this modification in Section  \ref{nonav}, in addition to given several technical lemmas.
The optimality  proofs for $m=6,5$ ($n=33, 22$) are then given in sections \ref{sec33}, \ref{sec22}, respectively.

\section{Continuously changing unavoidable configurations}
Optimality proofs for square packing utilize arguments based on resource starvation. Subsets of a containing square are associated with numerical resources  in such a way that each
packed box uses up a certain amount of resources (by intersecting the subset corresponding with the resource). The overall amount of resource available limits the number of boxes
that can be packed.

The proofs in \cite{goe} and many later publications are based on finite number of points, each of which has resource value $1$. In \cite{el}, resources were associated with line segments, such that the length of intersection between a box and the line segment determined the amount of resource allocated to the box. A more complex configuration in \cite{naga} uses a  combined system of (weighted) points, line segments, and a rectangular area.
 
 Our arguments will use a two-tier approach. We will first start out with systems of points containing too many resources for a direct proof. A new  technical 
 result (Theorem \ref{t:cont}) will allow us to use the flexibility in our initial systems to show that any potential packing must contain a local abundance of boxes. We will use
 this local ``over-concentration" of boxes to obtain a contradicting in combination with a second resource system based on a line segment.

We will start by stating several ``non-avoidance'' lemmas, which guarantee that a box will intersect particular type of subsets in its vicinity.     
\label{nonav}
\begin{lemma}[Friedman \cite{fried}, Stromquist \cite{strom2}] \label{triangle} Let $T$ be a triangle with sides of length at most 1. Then any box whose centre is in 
$T$ must contain one of the vertices of $T$.
\end{lemma}

\begin{lemma}[Friedman \cite{fried}, Stromquist \cite{strom2}] 
\label{rimevendistance}
Let $a \le 1$, $b\le 1$, and $a+2b \le 2\sqrt{2}$, then any box whose centre is in the rectangle 
$[0,a]\times[0,b]$ must intersect the $x$-axis, the point $(0,a)$ or the point $(a,b)$.
\end{lemma}
We will  use Lemma \ref{rimevendistance} in the cases of $a<2\sqrt{2}-2 \approx 0.828$, $b=1$ and 
$a=1$, $b< \sqrt{2}- \frac{1}{2} \approx 0.914$.

\begin{lemma}[Stromquist \cite{strom2}, \cite{strom}] 
\label{rimunevenpoints}
Let $2\sqrt{2}-2<a<1,\, 0<b<1$, and $(a,b)$ within a distance of $1$ from $(0,1)$. Moreover,
let $f(a)$ be the infimum of \begin{equation}\frac{\cos \theta}{1+\cos \theta}+\frac{1-a \cos \theta}{\sin \theta}\label{fa}
\end{equation} for 
$\theta \in (0,\frac{\pi}{4}]$. 
If $b<f(a)$, then any box whose centre is in the quadrilateral with vertices $(0,0),(0,1),(a,0),$ and $(a,b)$ must intersect
the $x$-axis, the point $(0,1)$ or the point $(a,b)$.
Moreover, the infimum of (\ref{fa}) is a minimum and is obtained at a value of $\theta$ satisfying
\begin{equation}\label{3rdpoly} 2 \cos^3 \theta -(2a+2) \cos^2 \theta + (a^2-2a+3) \cos \theta -(1-a^2)=0.\end{equation}
\end{lemma}

We will be using  Lemma \ref{rimunevenpoints} in the case $a=\frac{1}{2}\sqrt{3}, b=0.5$.

\begin{lemma}[Nagamochi \cite{naga}] If $l$ is a line that lies within a distance of $(\sqrt{2}-1)/2$ of the centre of a box $B$, then
$l$ will intersect $B$ with a length of more than $1$. \label{closeline}
\end{lemma}

\begin{lemma}[Stromquist \cite{strom2}] \label{l:paralines} Let $L_1$ and $L_2$ be two parallel lines of distance $d\le1$, and $B$ a box with 
its centre between them.  Then $B$ must intersect the 
two lines with a common length of intersection of at  least $\min\{1,2\sqrt{2}-2d\}$. 
\end{lemma}
The following lemma extends  Lemma \ref{rimunevenpoints} to values of $a$ smaller than $2\sqrt{2}-2$.
\begin{lemma} Let $0<a<2\sqrt{2}-2,\, 0<b\le1$, and $(a,b)$ within a distance of $1$ from $(0,1)$. 
Then any box whose centre is in the quadrilateral $Q$ with vertices $(0,0),(0,1),(a,0),$ and $(a,b)$ must intersect
the $x$-axis, the point $(0,1)$ or the point $(a,b)$.
\end{lemma}
{\bf Proof: } If $b\le \frac{1}{2}$ then the distance from $(0,0)$ to $(a,b)$ is less than $1$ and so the line segment between these points divides $Q$ into two triangles, all of whose sides have length at least one. The result now follows from Lemma \ref{triangle}.

So assume that $b >0.5$ and that the box $B$ does not intersect $(0,1)$ or the $x$-axis. Now the two line segments from $(0,0)$ to $(a,\frac{1}{2})$, and from 
$(0,1)$ to $(a,\frac{1}{2})$  divide $Q$ into three triangles, such that Lemma \ref{triangle} is applicable to each of them. By our assumption, the box $B$ must contain either $(a,b)$ or
 $(a,\frac{1}{2})$. In the first case, the lemma holds, so assume that $B$ contains  $(a,\frac{1}{2})$.

As the centre of $B$ is contained in $Q$, it is also contained in the larger rectangle $R$ with corners $(0,0)$, $(a,0)$, $(a,1)$, and $(0,1)$. Applying Lemma \ref{rimevendistance}
to $R$ yields that $B$ contains $(a,1)$. As $(a,b)$ lies on the line segment from $(a,\frac{1}{2})$ to $(a,1)$, it is contained in $B$. \hfill{$\Box$}\smallskip

We will use the lemma for $a=0.8$, $0.4 \le b \le 1$.
\begin{lemma} \label{l:pBl}
Let $l$ be a line and $P$ a point with a distance of  more than $0.51$ from $l$. If a box $B$ covers $P$ such that $P$ and the center of $B$ lie on opposite sides of $l$, 
than $B$ intersect $l$ with a length of intersection that exceeds $1$.  
\end{lemma}
{\bf Proof: }
The midpoint of B must lie within $0.505\sqrt{ 2}$ of $P$, and hence within a distance of $0.505 \sqrt{2}-0.51$ from the line $l$.
 As this value is less than $( \sqrt{2}-1)/2$, the results follows from Lemma 5.
\hfill{$\Box$}\smallskip

Consider a square $S$ of side length $l$ in the Euclidean plane (which in our cases we will take to be $[0,m]\times [0,m]$ for $m \in \{5,6\}$). A set of points $P \subset S$ is called \emph{unavoidable} if every box $B \subseteq S$ contains one of the points in $P$. In practice, we show unavoidability of $P$ by dividing $S$ into several regions $S_i$, so that by one of our unavoidability lemmas, any box with midpoint in $S_i$ must either intersect a point in $P$ or the boundary of $S$. If $S$ contains 
 an unavoidable set of $t$ points, it follows that no more than $t$ boxes can be packed into $S$, and hence $s(t+1)\ge l$.  
 
\begin{figure}[ht]
	\centering
		\includegraphics[trim=80 29 80 18,clip=true]{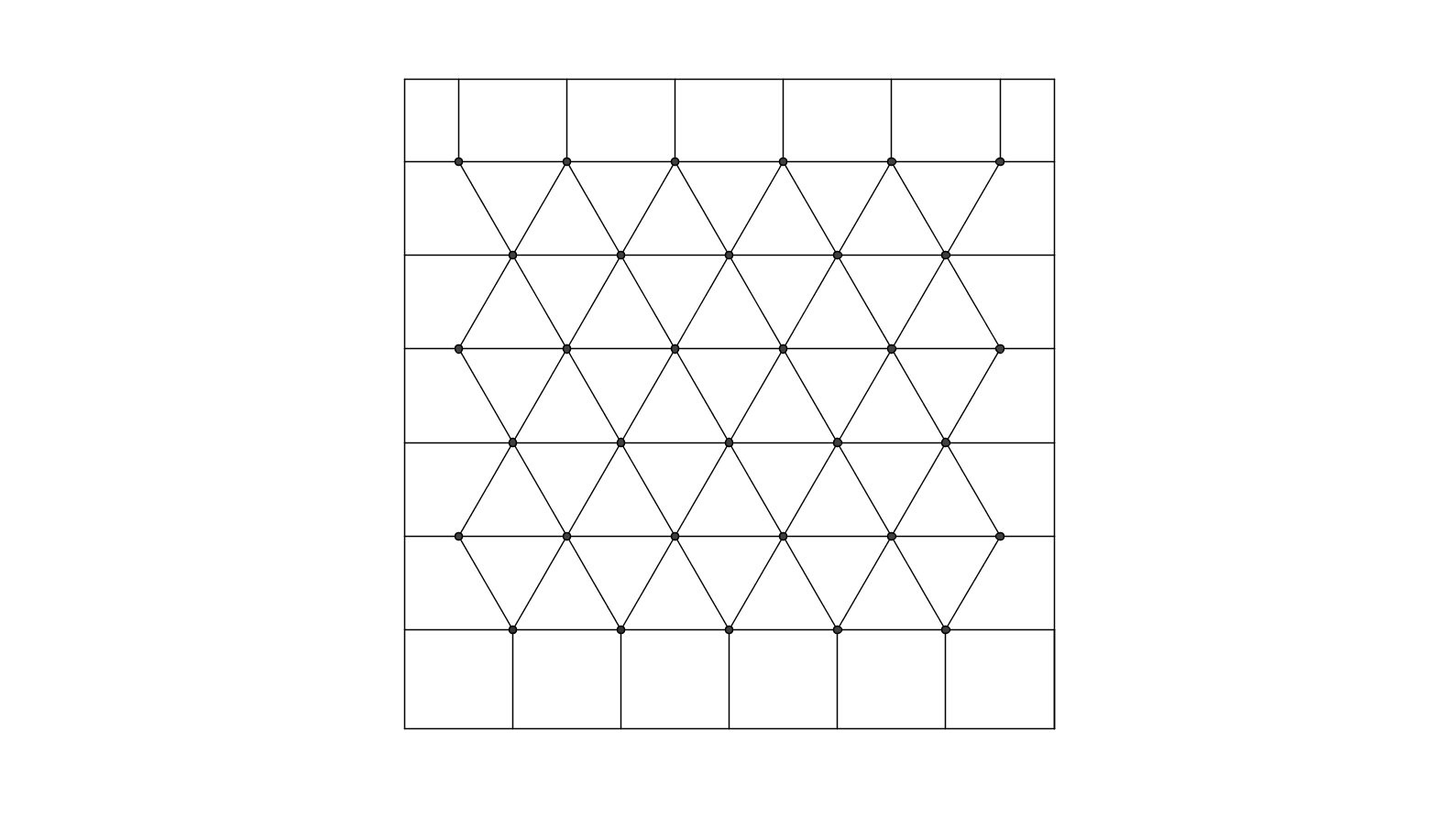}
	\caption{An unavoidable sets with 33 points}\label{6pat}
	
\end{figure}

Figure \ref{6pat} depicts an unavoidable set of points for the square $[0,6]\times [0,6]$. The points in the lowest row are 
$$\left(i,\sqrt{2}-\frac{1}{2}\right) \;\;\;\;\; i=1,2,\ldots,5,$$ 
and the remaining ones are 
arranged so that all shown triangles are equilateral of side length 1. 
Lemma \ref{triangle} is applicable to the triangles, Lemma \ref{rimevendistance} to the rectangles, and 
Lemma \ref{rimunevenpoints} to the remaining quadrilateral regions (with $a=\frac{\sqrt{3}}{2}$, 
$b=\frac{1}{2}$). 
 
 Figure \ref{6pat} is a variant of configurations used to show that $s(46)=7$ in \cite{b46} and to derive a lower bound on $s(11)$ in 
 \cite{strom}.  As the unavoidable set consists of $33$ points, it demonstrates the (known) result that $s(34)=6$. 
 
Note that the configuration in Figure \ref{6pat} contains a degree of flexibility. For example, we can obtain a different 
unavoidable configuration by deleting one of the points closest to the left hand side of the square and instead adding a different  point
a small amount further to the right. If there exists a packing of 
$33$ boxes, then each of them must contain exactly  one point in each configuration, and hence  one box must contain both the deleted and added point (and the line segment between them), an argument that has appeared in several previous proofs. The next theorem shows that this approach can be generalized to situations in which more than one point is moved at one time.

\begin{theo} \label{t:cont} Let $S$ be a square with a packing $\mathcal{P}$ of boxes, $I=[a,b]$, $t \in \mathbb{N}$, and $f_k:I \to S$ a collection of  continuous mappings, for $1\le k\le t$. Suppose further
that 
\begin{enumerate}\item for each $i \in I$, $F_i=\{f_k(i)|1\le k \le t\}$ is an unavoidable set of points;
\item if for some $1 \le k \le t$, $f_k(a)$ is not contained in a box of $\mathcal{P}$, then $f_k(i)=f_k(a)$ for all $i \in I$.
\item  if for some $1 \le k,l \le t$, $k\ne l$, $f_k(a)$ and $f_l(a)$ lie in the same box of $\mathcal{P}$, then $f_k(i)=f_k(a)$ for all $i \in I$.
\end{enumerate}
Then for all  $1 \le k \le t$, the image $f_k(I)$ will either lie entirely within one box, or completely outside any box. 
\end{theo}
{\bf Proof: }
If $f_k(a)$ is not contained in any box, then  $f_k$ is constant. Hence for the theorem to be wrong, there
must be  $1\le k\le t$, $i\in I$, such that $f_k(a)$ lies in some box $B_k$ while $f_k(i) \notin B_k$. Assume that this is indeed the case.

As boxes are open and $f_k$ is continuous, it
follows that there is a smallest such $i' \in I$ for which $f_k(i')$ lies outside $B_k$. Minimizing over all indices, we may assume w.l.o.g. that $i'$ is the smallest value of $i$ for which any $f_s(i)$ lies outside the box containing $f_s(a)$.   

Now, as $F_{i'}$ is an unavoidable set of points, there exist a $1 \le l\le k$, necessarily with $l \ne k$, such that $f_l(i') \in B_k$. Boxes are open, therefore
there exist an $\epsilon >0$ such that  $f_l(i'-\epsilon) \in B_k$. By the minimality of $i'$, it follows that $f_l(a) \in B_k$. However, now $f_l(a), f_k(a) \in B_k$, and so $f_k$ is constant by condition 3., contradicting that $f_k(i') \notin B_k$. The result follows.
\hfill{$\Box$}\smallskip

\section{The best possible packing of 33 unit squares}
\label{sec33}

\begin{theo}
\label{s33}
33 non-overlapping unit squares cannot be packed in a square of side length less than 6.
\end{theo}

\begin{figure}[ht]
	\centering
		\includegraphics[trim=80 29 80 20,clip=true]{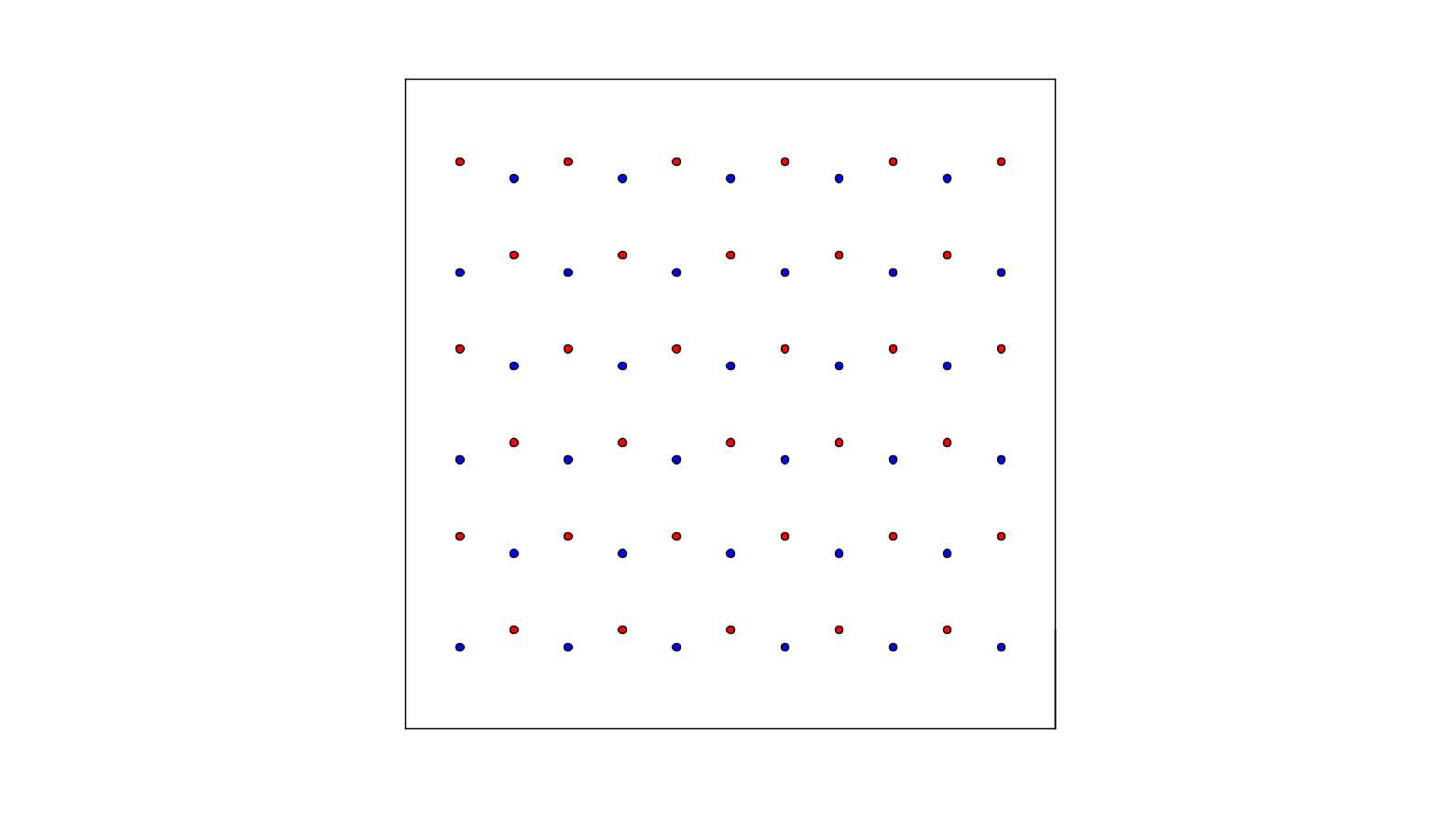}
	\caption{Two unavoidable sets with 33 points each}\label{fig:redblue}
\end{figure}

{\bf Proof: }
Let $S$ be the square $[0,6]^2$ and assume by way of contradiction that there is a packing $\mathcal{P}$ of 33 boxes into $S$.
Consider the collection of 33 red points and 33 blue points depicted in Figure 
\ref{fig:redblue}. The red points are the points from Figure \ref{fig:redblue}, while the blue points are obtained from the red ones by mirroring along the line $y=3$. It follows that both red and blue 
points form unavoidable sets of points. Hence each box in $\mathcal{P}$ will contain exactly one red and blue point.    

We will apply Theorem \ref{t:cont} twice. In the first instance we choose our values $f_k(i)$ so that $F_a$ is the set of red points from Figure \ref{fig:redblue}, while in the second case 
the start configuration will be the set of blue points. We will 
give an informal description of the other values of $f_i(k)$ by describing the ``movement" of the points configuration. Many of our movement will move an entire row of equally-colored points from Figure \ref{fig:redblue}. We will denote the red rows and blue rows by $r_1,\dots,r_6$, and $b_1,\dots,b_6$, respectively, where the rows are numbered from  bottom to top. 
For any two such rows $r,r'$, we denote by $v(r,r')$ their vertical distance.

We first note that as every red point and every blue points lies in exactly one box of $\mathcal{P}$, the second and third condition of Theorem \ref{t:cont} are automatically full-filled. 

 Now consider 
a simultaneous  vertical movement of a row $r_i$. Such a move will preserve the unavoidability of the   points configuration as long as, whenever they are defined, 
$v(r_i,r_{i-1}) \le \frac{1}{2}\sqrt{3}$, $v(r_i,r_{i+1})\le \frac{1}{2}\sqrt{3}$  and in addition, the vertical distance from $r_1$ and $r_6$  to the top or bottom edge of $S$, respectively, is 
at most  $\sqrt{2}-\frac{1}{2}$. 

With regard to a given configuration, for  $i=2,\dots ,5$, let $m_1=v(r_1,r_{2})$, $m_i=\max\{v(r_i,r_{i+1}), v(r_i,r_{i-1})\}$, for  $i=2,\dots ,5$, and $m_6=v(r_6,r_5)$.
Clearly, for $i=1,\dots,6 $, there exists a unique configuration $F_i$ that minimizes 
$m_i$ and that is reachable from $F_a$ by vertical movement of rows, such that unavoidability is preserved throughout.
Let 
$y_1, \dots, y_6$ be the second coordinate values of the points in the $i$-th row in $F_i$ (the exact values of $y_i$ can be easily calculated, but are not needed).  
As 
$$2\left(\sqrt{2}-\frac{1}{2}\right)+2\cdot  0.8 + 3\cdot\frac{1}{2}\sqrt{3}>6,$$ 
we note that in $F_i$, we have $v(r_i,r_{i-1}), v(r_i,r_{i+1})\le 0.8$, wherever  defined. 

Now consider $F_i$ for $i=2,4,6$. Here the the $i$-th row contains $6$ points, and  we may move the  $i$-th row horizontally to to the left and right, provided the distances
to all ``critical" points in the adjacent rows stays within $1$. As the adjacent rows have a vertical distance of less than $0.8$, it is easy to check that this allows for a movement of at least $0.1$ to either side. In addition, in this situation, we may move the left-most point of any such row horizontally to the right, until it reaches the point $(1,y_i)$. This maximal reflection is possible, as the distance to the adjacent row is less than  $2\sqrt{2}-2\ge 0.8$.

We now proceed as follows. We move the point configuration to one of the $F_i$. If the $i$-th row contains $5$ points, we note that one of the points occupies the point $(1,y_i)$. 
If the $i$-th row contains $6$ points, we move the row $0.1$ to the left and back to the right. Finally, move the left-most point of row $i$ from $(0.5,y_i)$ to $(1,y_i)$ and back, and note that this point, combining both movements, has moved over the line segment from $(0.4,y_i)$ to $(1,y_i)$. We repeat this procedure for all $F_i$.

By Theorem \ref{t:cont}, the line segments $[0.4,1]\times \{y_i\}$ for $i\in \{2,4,6\}$ lie within the same box of $\mathcal{P}$. Moreover, as these sets 
are the result of movement from  different points of the base configuration, they, as well as the point sets $\{(1,y_i)\}$ from the movement of the $5$-point rows lie in different boxes.

We now repeat this procedure for the blue points, noting that we obtain the same values $y_i$ as in the case of the red points. The blue points will have a 
$6$-point row where the red points have a $5$-point row and vice versa. We can conclude that for $i\in\{1,3,5\} $, the line segment $[0.4,1]\times \{y_i\}$ lies within one box of $\mathcal{P}$. 

Hence, taking both configurations together, the segments $[0.4,1]\times \{y_i\}$ lie each in one box of $\mathcal{P}$ for $i=1,\dots, 6$. Moreover, these segments all lie in different boxes 
as the points $(1,y_i)$ are all within the movement of different points from the red base configuration. 

It follows that in $\mathcal{P}$, there are $6$ distinct boxes $B_1,\dots,B_6$ such that  $B_i$ covers $(0.4,y_i)$. Now let $l$ be the line segment from
$(\sqrt{2}-\frac{1}{2},0)$ to $(\sqrt{2}-\frac{1}{2},6)$. We are interested in the length of the intersection of $B_i$ and $l$. If  the midpoint of $B_i$ lies on the same side of $l$ as $(0.4, y_i)$, then 
this intersection exceeds $1$ by Lemma \ref{l:paralines} (with $d=\sqrt{2}-\frac{1}{2}$). If the midpoint of $B_i$ lies on $l$ or on the side of $l$ opposite from $(0.4,y_i)$, then the intersection has length larger 
than $1$ by Lemma  \ref{l:pBl}. Hence all six boxes intersect $l$ with a length larger than $1$. However, the length of $l$ is $6$, for a contradiction.

Hence $33$ boxes cannot be packed in a square of side length $6$ and so $s(33)\ge 6$.\hfill{$\Box$}\smallskip

\begin{cor}The trivial packing is optimal for packing   $33$ unit squares in a square, and we have that  $s(33)=6$. 
\end{cor}

\section{The best possible packing of 22 unit squares}
\label{sec22}
\begin{theo}
\label{th22}
$22$ non-overlapping unit squares cannot be packed in a square of side length less than $5$.
\end{theo}

\begin{figure}[h]
	\centering
		\includegraphics[scale= 0.33,trim=180 44 180 27,clip=true]{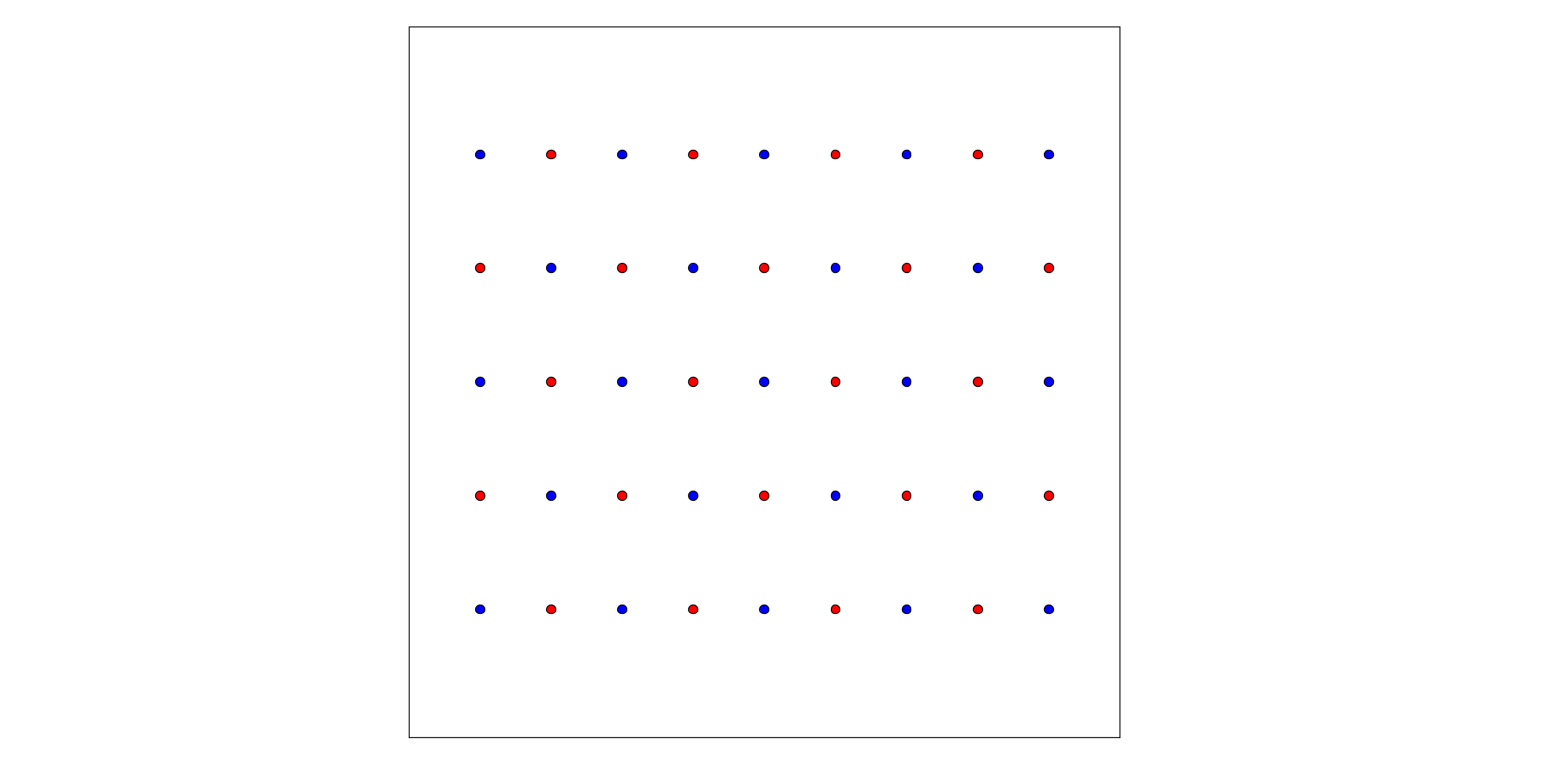}
		\caption{Unavoidable configurations of 22 red and 23 blue points}\label{fig:bluered22}
\end{figure}

{\bf Proof: } Let $S$ be the square $[0,5]^2$, and assume, by way of contradiction, that there exist a packing $\mathcal{P}$ of $22$ boxes into $S$.
Consider the configuration of $22$ red and $23$ blue points  shown in Figure \ref{fig:bluered22}. Together, the points of both colours are exactly the elements the set 
$\{0.5,1,1.5,2,2.5,3,3.5,4,4.5\} \times \{0.9,1.7,2.5,3.3,4.1\}$. We will denote  the rows of red and blue points by $r_1,\dots,r_5$, and $b_1,\dots, b_5$, with numbering from bottom
to top, and we let $y_1,\dots,y_5$ be the value of their second coordinate.

By the same basic argument as applied to  Figure \ref{fig:redblue} we can show that both the set of red points and the set of blue points form unavoidable sets. It follows that each red points lies 
in exactly one box from $\mathcal{P}$. In the case of the blue points, either exactly 
one blue point is not contained in a box of $\mathcal{P}$ or exactly one box of $\mathcal{P}$ contain two blue points. In the latter case, the two blue points that lie in the same box must be within a distance of each other that is smaller than the diagonal of a maximal size box, which is  $1.01 \sqrt{2}$.
Thus, by the symmetry in our blue point configuration, we may assume that all blue points with first coordinate value
 lower than $2$ lie in a box of $\mathcal{P}$ that does not contain any other blue points.  

As in the proof of Theorem \ref{s33}, we will apply Theorem \ref{t:cont} twice,  with the red and blue points as the respective starting configurations. Once again we 
describe the function $f_k$ informally in terms of the movement of points. In the case of the red points, we will first move $r_2$ and $r_4$ (i.e. the red rows 
containing $5$ points) horizontally a distance of $0.1$
 to the left and back, followed by 
a movement of the left-most point in each such row horizontally to the right until it reaches a first coordinate value of  $1$.  These movement preserve
 unavoidablility, as the vertical distance between adjacent rows is $0.8$ in all cases. As in Theorem \ref{s33} we record that different boxes of $\mathcal{P}$ contain the line segments $[0.4,1]\times \{1.7\}$ and $[0.4,1]\times \{3.3\}$, and that these boxes must also be different from the boxes containing the (stationary) points $(1,0.9)$, $(1,2.5)$, and $ (1,4.1)$. 
 
For  the blue points, we want repeat these movement with rows $b_1,b_3,$ and $b_5$. However, we need to modify our procedure, as Theorem \ref{t:cont} requires the a blue point remains stationary 
if it does not lie in a box of $\mathcal{P}$ or is in a box of $\mathcal{P}$ that contains a different blue point. If row $r_i$, $i\in \{1,3,5\}$ does not contain such point, we move it horizontally a distance of 
$0.1$ and back to its original position. We then move its leftmost point horizontally from $(0.5,y_i)$ to $(1,y_i)$. In case that $r_i$ does contain such point, we just move its leftmost point from $(0.5,y_i)$ to $(1,y_i)$. We note that the last case can only happen for one of the rows $b_1,b_3,b_5$, because if there are two such exceptional blue points, they must lie in the same box, and hence (due to the size of the boxes) in either the same or adjacent rows.  

By Theorem \ref{t:cont} the trajectory of each of the leftmost points of $b_1,r_2,b_3,r_4,b_5$ lies completely within a box of $\mathcal{P}$. Moreover each of these trajectories intersects a different red point from the initial configuration, and hence the trajectories lie within different boxes. 

We can conclude that there are $5$ boxes $B_1, \dots, B_5$ in $\mathcal{P}$
such that the line segments $[0.4,1] \times \{y_i\}$ lie completely within $B_i$, except that at most one of $B_1,B_3,B_5$ might only cover the line segment  
$[0.5,1] \times \{y_i\}$.  

\begin{figure}[h]
	\centering
		\includegraphics[scale= 0.6,trim=0 0 70 38,clip=true]{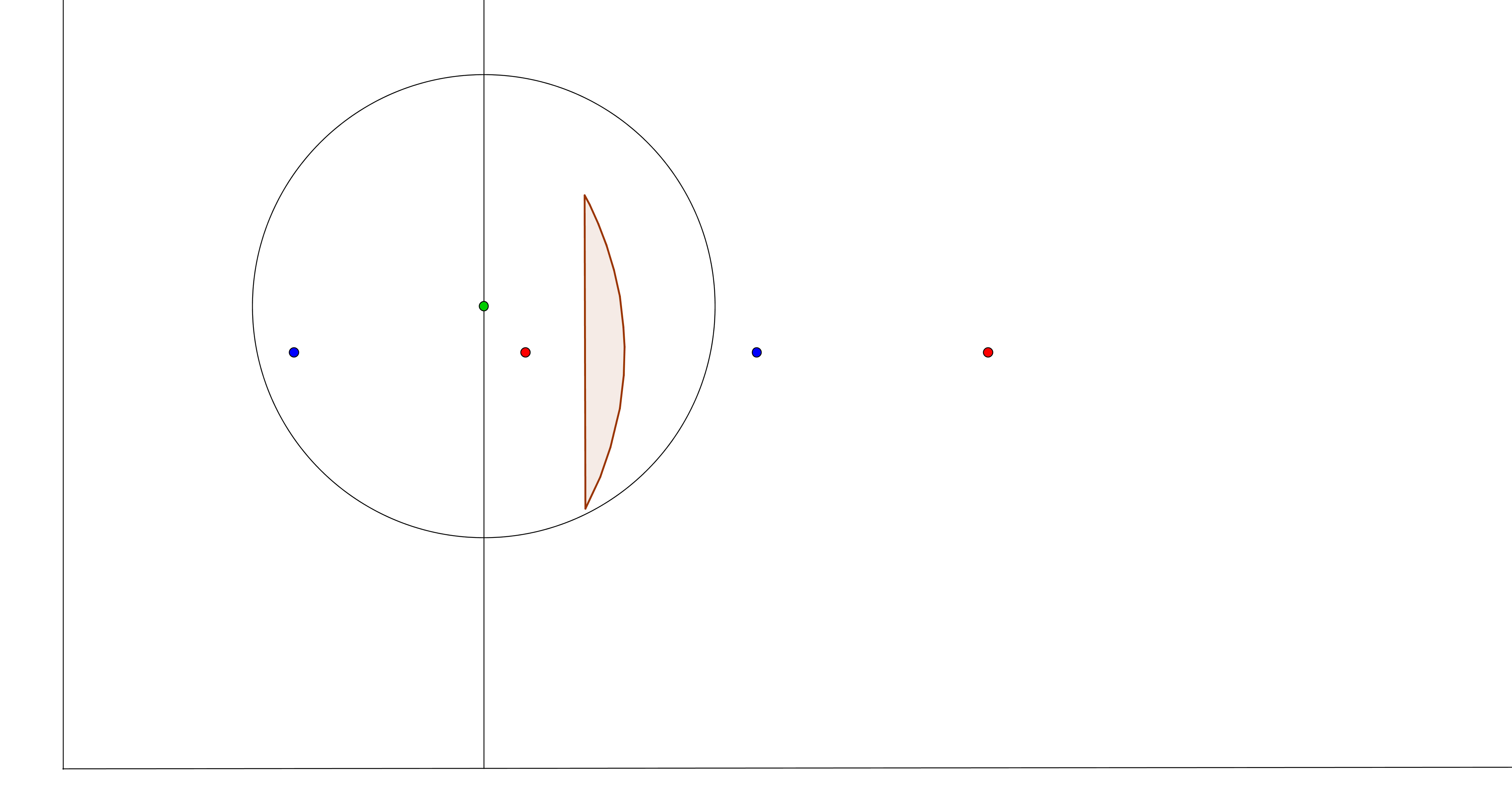}
		\caption{The midpoint of the box $B_1$ must lie in the shaded area}\label{fig:lowline}
\end{figure}

Let $l$ be the line segment $\{\sqrt{2}-\frac{1}{2}\} \times [0,5]$. As in Theorem \ref{s33} we can check that if $B_i$ contains $[0.4,1]\times \{y_i\}$, it  intersect
$l$ with a length of intersection  that exceeds $1$. As $l$ has length $5$, one of $B_1,B_3,B_5$ does not  cover the entire  line segment from $[0.4,1]\times\{y_i\}$. 
We consider $2$ cases:
\begin{enumerate}
\item First assume that $B_3$ does not completely cover $[0,4]\times\{y_3\}$. As there was at most one exceptional row, $B_1,B_2,B_4,$ and $B_5$ all intersect $l$ with a length of intersection 
exceeding $1$. It follows that $B_3 \cap l \subseteq \{\sqrt{2}-\frac{1}{2}\}\times (2,3)$, and so $B_3$ does not cover the points $(\sqrt{2}-\frac{1}{2},2)$ and $ (\sqrt{2}-\frac{1}{2},3)$. In Figure \ref{fig:midline}, these two points are depicted in green.

Let $m$ be the midpoint of $B_3$. The location of $m$ is constraint as follows: $m$ must lie on the right side of $l$  and separated from it by  a distance of at least $\frac{1}{2}\sqrt{2}-\frac{1}{2}$, for
otherwise the length of intersection of $B_3$ and $l$ would  exceed $1$ by either Lemma \ref{l:paralines} or Lemma \ref{closeline}. As $B_3$ does not cover $(\sqrt{2}-\frac{1}{2},2)$ or $(\sqrt{2}-\frac{1}{2},3)$, $m$ cannot be 
within a distance of $0.5$ from either of these points. Finally, as $B_3$ covers $[0.5,1] \times \{y_3\}$, the distance from $m$ to $(0.5,y_3)$ must be smaller than half the 
diagonal of a maximal box, i.e. smaller than $0.505\sqrt{2}$. Figure \ref{fig:midline} shows the remaining possible locations of $m$ as a shaded area. The area is bounded 
by line and circle segments that intersect in 
$4$ points with approximate coordinates $(1.12, 2.5 \pm 0.05), \,    (1.2,2.5\pm 0.1)$.

An easy calculation shows that
the entire area is within a distance of $0.5$ from the point $(1.5,y_3)=(1.5,2.5)$. In Figure \ref{fig:midline}, this distance is indicated by a circle. It follows that $B_3$ also covers the point $(1.5, 2.5)$.
Hence in our initial configuration of points, the box $B_3$ covers two blue points, namely those at $(0.5,2.5)$ and $(1.5,2.5)$. However, this contradict our assumption that all blue points with a second coordinate value smaller than $2$ do not share a box with another blue point.

\begin{figure}[h]
	\centering
		\includegraphics[scale= 0.6,trim=0 40 70 40,clip=true]{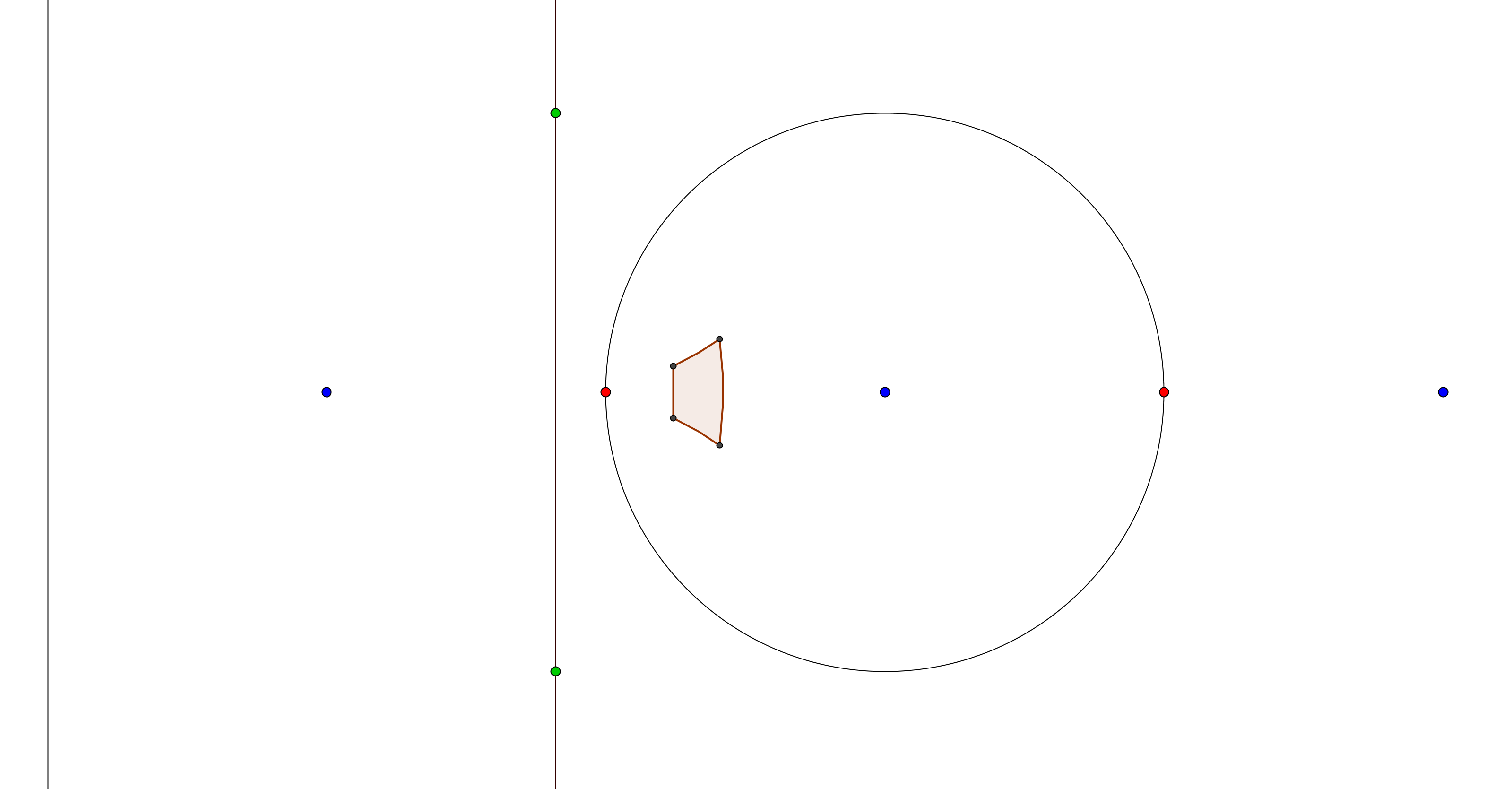}
		\caption{The midpoint of the box $B_3$ must lie in the shaded area}\label{fig:midline}
\end{figure}
\item Assume that for one $i \in \{1,5\}$, $B_i$ does not cover the entire line segment $[0.4,1] \times \{y_i\}$. By symmetry, we may assume that this is the case  for $i=1$. Then 
$B_i$ covers $[0.4,1]\times \{y_i\}$ for $i=2,\dots, 5$, and, as in the previous case, this implies that each such $B_i$ intersects the line segment $l$ with a length of intersection larger than 
$1$. It follows that the point $(\sqrt{2}-\frac{1}{2},1)$ is denied to $B_1$. This point is depicted green in Figure \ref{fig:lowline}. 

Let $m$ be the midpoint of $B_1$. As before we can conclude that $m$ lies on the opposite side of $l$ from the point $(0.5,\sqrt{2}-\frac{1}{2})$, with a distance of at least $\frac{1}{2}\sqrt{2}-\frac{1}{2}$ from $l$, but within a distance of $0.505\sqrt{2}$ of $(\frac{1}{2},\sqrt{2}-\frac{1}{2})$.  These constraints intersect  at approximately $(1.13, 0.56)$ and   
$ (1.13,1.24)$.

The resulting area is depicted in Figure \ref{fig:lowline} and lies completely within a distance of $\frac{1}{2}$ from the point $(\sqrt{2}-\frac{1}{2},1)$. It follows that $(\sqrt{2}-\frac{1}{2},1)$ lies in 
$B_1$. However, the point is denied to $B_1$, for a contradiction.
\end{enumerate} 
In either case we get a contradiction. It follows that the packing  $\mathcal{P}$ does not exists, and hence $s(22)\ge 5$.\hfill{$\Box$}\smallskip

\begin{cor}The trivial packing is optimal for packing   $22$  unit squares in a square, and we have that $s(22)=5$. 
\end{cor}


\end{document}